\documentclass[11pt]{article}

\usepackage{times}
\usepackage{alltt}
\usepackage{amssymb}

\usepackage{amssymb,amsmath,latexsym}    

\usepackage{array}      

\newtheorem{thm}{Theorem}

\newtheorem{lem}{Lemma}

\newtheorem{defn}{Definition}
\newtheorem{exmp}{Example}
\newtheorem{rem}{Remark}

\title{Quadratic Function Fields with Exponent Two Ideal Class Group}
\author{Victor Bautista-Ancona, Javier Diaz-Vargas \footnote {Partially supported by Programa de Impulso y Orientaci\'{o}n a la Investigaci\'{o}n, Universidad Aut\'{o}noma de Yucat\'{a}n}}
\date{Facultad de Matem\'{a}ticas. 
      Universidad Aut\'{o}noma de Yucat\'{a}n.}

\begin{document}

\maketitle

\section*{Abstract}
It has been shown by Madden that there are only finitely many quadratic extensions of $k(x)$, $k$ a finite field, in which the ideal class group has exponent two and the infinity place of $k(x)$ ramifies. We give a characterization of such fields that allow us to find a full list of all such field extensions for future reference. In doing so we correct some errors in earlier published literature.

\section{Introduction}

The study of function fields over finite constant fields with exponent two ideal class group is made by considering those fields with exponent two group of divisor classes of degree zero. If $K$ is a quadratic extension of $k(x)$, in which the infinite place ramifies, then the ideal class group and the group of divisor classes of degree zero are identical.  Notice that any place of degree $1$ may act as the infinite place. Also, exponent $2$ class group means class number not one, so genus not zero. Using these facts, Madden \cite{qff} has proved:

\begin{thm}[Madden]
\label{madden}If $K$ is a quadratic extension of $k(x)$ of genus $g$, where $k$ is a finite field of order $q$, in which the infinite place ramifies, and if the ideal class group of $K$ has exponent two, then:\newline
$\left( i\right) $ $q=9$, $g=1$;\newline
$\left( ii\right) $ $q=7$, $g=1$;\newline
$\left( iii\right) $ $q=5$, $1 \leq g \leq 2$;\newline
$\left( iv\right) $ $q=4$, $1 \leq g \leq 2$;\newline
$\left( v\right) $ $q=3$, $1 \leq g \leq 4$;\newline
$\left( vi\right) $ $q=2$, $1 \leq g \leq 8$.
\end{thm}

This theorem is the basis for our search of all the quadratic function fields with exponent two ideal class group.

The content of the paper is the following. The main result that relates the class number with the number of places that ramify is given in the second section. Based on this fact, we conclude that in an exponent two quadratic function field, where the infinite place ramifies, the class number, a power of 2, can not be greater than $2^{5}$. Some facts needed about Artin-Schreier extensions, elliptic and hyperelliptic function fields are established in the section three. The next section introduces some more specific topics and at the end of the section we give, for some genus, a formula for the class number in terms of the number of places up to degree equal to the genus. The section 6 gives some results about function fields with class number 2, 4, 8, 16 and 32. There, we give necessary conditions over the finite fields and the genus of the function fields to have the desirable class number. Section 7 contains a brief discussion of how to find the examples that we are looking for, and finally, we give the full list of function fields whose ideal class group has exponent two, up to $\mathbb{F}_{q}$-isomorphism.

\section{Characterization of exponent two quadratic extensions}

\begin{thm}
\label{caracterizacion}Let $K$ be a quadratic extension of $k(x)$, where $k$ is a finite field of order $q$, in which the infinite place ramifies, and $h$ be its class number. Then, the group of divisor classes of degree zero has exponent two if and only if $h=2^{t-2} $ or $h=2^{t-1}$ depending on whether $q \equiv 1 \pmod{2}$ or not ($t$ is the number of places of $k(x)$ that ramify in $K$).
\end{thm}

{\bf Proof. }
In a quadratic extension $K$ of $k(x)$, in which the group of divisor classes of degree zero has exponent two, all the classes in this group are ambiguous, i.e., they are fixed by the Galois group of $K/k(x)$. Thus the class number $h_{K}$ of $K$ is equal to the ambiguous class number $\left| J_{K}^{G}\right| $. Since the infinite place ramifies,
because of theorem $9$ in (\cite{adc}, page 167), we have that: 
\begin{eqnarray*}
q & \equiv &0 \pmod{2} \Rightarrow \left| J_{K}^{G}\right| =h_{k\left(x\right) }2^{t-1}, \\
q & \equiv &1 \pmod{2} \Rightarrow \left| J_{K}^{G}\right| =h_{k\left(x\right) }2^{t-2}.
\end{eqnarray*}
From the fact that $\left| J_{K}^{G}\right| =h_{K}$ and $h_{k(x)}=1$ the result follows.
The converse follows easily from theorem $12$ in (\cite{adc}, page 169). $\blacksquare $

\section{Some background material}

When the characteristic is two, we use some properties of \emph{Artin-Schreir extensions} of $K/\mathbb{F}_{q}(x)$ and the \emph{Hasse Normal Form} associated with them. Hasse gave much information about the arithmetic of such extensions. The Hasse Normal Form will help us to find equations for the quadratic extensions of $K/\mathbb{F}_{q}(x)$. When the characteristic is different from two, we use \emph{elliptic and hyperelliptic function fields }according to if the genus is one or greater or equal to two, respectively. For proofs, see (\cite{aff}, pages 115-118, 186-195).

\subsection{Artin-Schreier extensions}

\begin{thm}
\label{as}Let $F/k$ be a function field of characteristic $p>0$. Suppose that $u\in F$ is an element which satisfies the following condition: 
\begin{equation}
u\neq w^{p}-w \textmd { for all } w\in F.  \label{artin cond}
\end{equation}
Let 
\[K=F(y) \textmd { with } y^{p}-y=u.\]
Such an extension $K/F$ is called an \emph{Artin-Schreier extension} of $F$. For $\mathcal{P}\in \mathbb{P}_{F}$ (the set of places of $F$), we define the integer $m_{\mathcal{P}}$ by 
\[m_{\mathcal{P}}:=\left\{
\begin{array}{cc}
m & 
\begin{array}{c}
\textmd { if there is an element } z\in F \textmd { satisfying } \\ 
v_{\mathcal{P}}\left( u-\left( z^{p}-z\right) \right) =-m<0 \textmd { and } m \not{\equiv} 0 \pmod{p},
\end{array}
\\ 
-1 & \textmd { if }v_{\mathcal{P}}\left( u-\left( z^{p}-z\right) \right) \geq 0
\textmd { for some }z\in F.
\end{array}
\right. 
\]
We then have:\newline
(a) $\mathcal{P}$ is unramified in $K/F$ iff $m_{\mathcal{P}}=-1$.\newline
(b) $\mathcal{P}$ is totally ramified in $K/F$ iff $m_{\mathcal{P}}>0$.%
\newline
(c) If at least one place $\mathcal{Q}\in \mathbb{P}_{F}$ satisfies $m_{%
\mathcal{Q}}>0$, then $k$ is algebraically closed in $K$, and 
\[
g_{K}=p\cdot g_{F}+\frac{p-1}{2}\left( -2+\sum_{\mathcal{P}\in \mathbb{P}%
_{F}}\left( m_{\mathcal{P}}+1\right) \cdot \deg \mathcal{P}\right) , 
\]
where $g_{K}$ (resp. $g_{F})$ is the genus of $K/k$ (resp. $F/k$).
\end{thm}

\begin{rem}
\label{cond suf as}Suppose that there exists a place $\mathcal{Q}\in \mathbb{P}
_{F}$ with 
\[
v_{\mathcal{Q}}\left( u\right) <0\textmd { and }v_{\mathcal{Q}}\left( u\right) \not{\equiv}0 \pmod{p}. 
\]
Then $u$ satisfies condition $(\ref{artin cond})$ of the previous theorem.
\end{rem}

\subsection{The Hasse Normal Form}

Furthermore, if $F=k(x)$ the rational function field, the partial fraction decomposition of $u$ allows us to adjust the function $y$, so that $u$ is in the Hasse Normal Form; that is, $u$ satisfies the following conditions: 
\[u=\frac{q\left( x\right) }{\prod_{i}p_{i}\left( x\right) ^{\gamma _{i}}},\]
where
\begin{enumerate}
\item  $p_{i}\left( x\right) $ are irreducible polynomials of $k(x)$.
\item  $\gamma _{i}$ are positive integers relatively prime to $p$.
\item  $q(x)$ is relatively prime to the denominator, and
\item  $\deg u=\deg q\left( x\right) -\deg \left( \prod_{i}p_{i}\left(x\right) ^{\gamma _{i}}\right) $ is either negative, zero, or relatively prime to $p$.
\end{enumerate}

Then, if $\mathcal{P}_{i}=\mathcal{P}_{p_{i}\left( x\right) }$ is the place of $k\left( x\right) $ associated to $p_{i}\left( x\right) $, and $\mathcal{P}=\mathcal{P}_{\infty }$ is the infinite place of $k\left( x\right) $,
\begin{enumerate}
\item  $v_{\mathcal{P}_{i}}\left( u\right) =-\gamma _{i}$,
\item  $v_{\mathcal{P}}\left( u\right) =-\deg u$.
\end{enumerate}

So, if $\deg u>0$ and relatively prime to $p$, because of the remark $\ref{cond suf as}$, we can apply the theorem $\ref{as}$ to see that the only places that ramify in $K=k\left( x,y\right) $ are $\mathcal{P}_{i}$ and $\mathcal{P}$. In addition, if $k=\mathbb{F}_{2}$, so that $p=2$ and we are looking at curves with 
\begin{eqnarray*}
m_{\mathcal{P}_{i}} &=&\gamma _{i}, \\
m_{\mathcal{P}} &=&\deg u,
\end{eqnarray*}
the genus formula in the theorem $\ref{as}$ specializes to: 
\[g=-1+\frac{1}{2}\left( \sum_{i}\left( \gamma _{i}+1\right) \cdot \deg p_{i}\left( x\right) +\left( \deg u+1\right) \right)\]
where $g=g_{K}$ and $\gamma _{i}$, $\deg u$ are odd positive integers.

Assume now that $g>0$. Then, 
\[2\left( g+1\right) =\sum_{i}\left( \gamma _{i}+1\right) \cdot \deg p_{i}\left( x\right) +\left( \deg u+1\right) .\]
But $\deg u=\deg q\left( x\right) -\sum_{i}\gamma _{i}\deg p_{i}\left(x\right) $ implies that 
\begin{eqnarray*}
2\left( g+1\right) &=&\sum_{i}\deg p_{i}\left( x\right) +\deg q\left(x\right) +1\textmd {, i.e., } \\
\deg q\left( x\right) &=&2g+1-\sum_{i}\deg p_{i}\left( x\right) .
\end{eqnarray*}
But if we want $\deg u>0$, that implies that 
\begin{eqnarray}
\deg q\left( x\right) &>&\sum_{i}\gamma _{i}\cdot \deg p_{i}\left( x\right) 
\textmd {, i.e., }  \nonumber \\
2g+1 &>&\sum_{i}\left( \gamma _{i}+1\right) \cdot \deg p_{i}\left( x\right) 
\textmd {. }  \label{mayor 0}
\end{eqnarray}
If $\deg p_{i}\left( x\right) =n$ and $\gamma _{i}=1$ but $\gamma _{j}=0$ for $j\neq i$, then $2g+1>2n$. Therefore, $n<\frac{2g+1}{2}$ or $n\leq g$. Also notice that if $\deg p_{i}\left( x\right) =g$, then if $\gamma _{i}\neq 0$ and $\gamma _{j}\neq 0$ for $j\neq i$, 
\[\left( \gamma _{i}+1\right) g+\left( \gamma _{j}+1\right) \deg p_{j}\geq 2g+2.\]
Therefore, if $\gamma _{i}\neq 0$, $\gamma _{j}=0$ for all $j\neq i$. Similarly, if $g>0$, $\gamma _{i}=1$.

\subsection{Elliptic function fields}

\begin{defn}
A function field $K/k$ (where $k$ is the full constant field of $K$) is said
to be an \emph{elliptic function field} if the following conditions hold:%
\newline
(1) the genus of $K/k$ is $g=1$, and \newline
(2) there exists a divisor $A$ with $\deg A=1$.
\end{defn}

\begin{thm}
\label{nec eff}Let $K/k$ be an elliptic function field. If char $k\neq 2$,there exist $x$, $y\in K$ such that $K=k\left( x,y\right) $ and 
\[y^{2}=f\left( x\right) \in k\left[ x\right]\]
with a square-free polynomial $f\left( x\right) \in k\left[ x\right] $ of degree $3$.
\end{thm}

Now, an ``algorithm'' for to find our examples is the following.

\begin{thm}
\label{suf eff}Let char $k\neq 2$. Suppose that $K=k\left( x, y\right) $ with 
\[ y^{2}=f\left( x\right) \in k\left[ x\right] , \]
where $f\left( x\right) $ is a square-free polynomial of degree $3$. Consider the decomposition $f\left( x\right) =c\prod_{i=1}^{r}p_{i}\left(x\right) $ of $f\left( x\right) $ into irreducible monic polynomials $p_{i}\left( x\right) \in k\left[ x\right] $ with $0\neq c\in k$. Denote by $\mathcal{P}_{i}\in \mathbb{P}_{k\left( x\right) }$ the place of $k\left(x\right) $ corresponding to $p_{i}\left( x\right) $, and by $\mathcal{P}_{\infty }\in \mathbb{P}_{k\left( x\right) }$ the pole of $x$. Then the following holds:\newline
($1$) $k$ is the full constant field of $K$, and $K/k$ is an elliptic function field.\newline
($2$) The extension $K/k\left( x\right) $ is cyclic of degree $2$. The places $\mathcal{P}_{1}$, $\ldots $, $\mathcal{P}_{r}$ and $\mathcal{P}_{\infty }$ are the only ones which are ramified in $K/k\left( x\right) $; each of them has exactly one extension in $K$, say $\mathcal{Q}_{1}$, $\ldots $, $\mathcal{Q}_{r}$ and $\mathcal{Q}_{\infty }$, and we have that the indices of ramification are $e\left( \mathcal{Q}_{j}\mid \mathcal{P}_{j}\right) =e\left( \mathcal{Q}_{\infty }\mid \mathcal{P}_{\infty }\right) =2$, $\deg \mathcal{Q}_{j}=\deg \mathcal{P}_{j}$ and $\deg \mathcal{Q}_{\infty }=1$.
\end{thm}

\subsection{Hyperelliptic function fields}

\begin{defn}
A \emph{hyperelliptic function field} over $k$ is a function field $K/k$ of genus $g\geq 2$ which contains a rational subfield $k\left( x\right) \subset K$ with $\left[ K:k\left( x\right) \right] =2$.
\end{defn}

\begin{lem}
(a) A function field $K/k$ of genus $g\geq 2$ is hyperelliptic if and only if there exists a divisor $A$ with $\deg A=2$ and $\dim A\geq 2$.\newline
(b) Any function field $K/k$ of genus $2$ is hyperelliptic.
\end{lem}

We can provide an explicit description of $K/k$.

\begin{thm}
\label{hyper eff}Assume that char $k\neq 2$. \newline
(a) Let $K/k$ be a hyperelliptic function field of genus $g$. Then there exist $x$, $y\in K$ such that $K=k\left( x,y\right) $ and 
\begin{equation}
y^{2}=f\left( x\right) \in k\left[ x\right]  \label{ecua}
\end{equation}
with a square-free polynomial $f\left( x\right) $ of degree $2g+1$ or $2g+2$.\newline
(b) Conversely, if $K=k\left( x,y\right) $ and $y^{2}=f\left(x\right) \in k\left[ x\right] $ with a square-free polynomial $f\left(x\right) $ of degree $m>4$, then $K/k$ is hyperelliptic of genus 
\[g=\left\{ 
\begin{array}{lll}
\frac{m-1}{2} &  & \textmd {if }m\equiv 1 \pmod{2}, \\ 
\frac{m-2}{2} &  & \textmd {if }m\equiv 0 \pmod{2}.
\end{array}
\right. 
\]
(c) Let $K=k\left( x,y\right) $ with $y^{2}=f\left( x\right) $ as in $(\ref{ecua})$. Then the places $\mathcal{P}\in \mathbb{P}_{k\left( x\right) }$ which ramify in $K/k\left( x\right) $ are the following: 
\begin{eqnarray*}
\textmd {all zeros of }f \left( x\right) \textmd { if }\deg f\left( x\right)&\equiv &0 \pmod{2}, \\
\textmd {all zeros of }f \left( x\right) \textmd { and the pole of }x\textmd { if } \deg f\left( x\right) &\equiv &1\pmod{2}.
\end{eqnarray*}
Hence, if $f\left( x\right) $ decomposes into linear factors, exactly $2g+2$ places of $k\left( x\right) $ are ramified in $K/k\left( x\right) $.\label{hyperelliptic}
\end{thm}

\section{ Upper bound}

From the theorem $\ref{caracterizacion}$ and the discussion about the Hasse Normal Form we get the following result.

\begin{thm}
\label{menorque}Let $K$ be a quadratic extension of $k(x)$, where $k$ is a finite field of order $q$, in which the infinite place ramifies and whose group of divisor classes of degree zero has exponent two, and let $h$ be its class number, a power of $2$. Then $h\leq 2^5$.
\end{thm}

{\bf Proof. }
 If $q=2$, suppose that $h=32$. From the theorem $\ref{caracterizacion}$, $h=2^{t-1}$; hence, $t=6$, i.e., there exist $6$ places of $\mathbb{F}_{2}\left(x\right) $ that ramify in $K$. Since the infinite place ramifies, then $5$ other places different from the infinite place must ramify. Recall that there are two irreducible polynomials of degree $1$ over $\mathbb{F}_{2}$, namely, $x$ and $x+1$; one irreducible polynomial of degree $2$, $x^{2}+x+1$; two polynomials of degree $3$, $x^{3}+x^{2}+1$ and $x^{3}+x+1$. Because of the equation $(\ref{mayor 0})$, taking $\gamma _{i}=1$, we have that $2g+1>2\left(1+1+2+3+3\right) =20$. From here, $g\geq 10$ (a contradiction with theorem $\ref{madden}$). Hence, for $q=2$, $h\leq 16$.

If $q=3$, suppose that $h=64$. From the theorem $2$, $h=2^{t-2}$. From here, $t=8$. So, $7$ places different from the infinite place must ramify. There exist three irreducible polynomials over $\mathbb{F}_{3}$ of degree $1$, namely, $x$, $x+1$ and $x+2$; three irreducibles of degree $2$, $x^{2}+1$, $x^{2}+x+2$, $x^{2}+2x+2$. Since the characteristic $\neq 2$, by theorem $\ref {hyperelliptic}$, part (c), $2g+1=1+1+1+2+2+2+4=13$, i.e., $g=6$ (again a contradiction with theorem $\ref{madden}$). From here, for $q=3$, $h\leq 32$.

If $q=4$, we know that $h\leq 2^{g}$, but by theorem $\ref{madden}$, $h\leq 4 $.

If $q=5$, from $h\leq 2^{2g}$ and theorem $\ref{madden}$, we conclude that $h\leq 16$.

Finally, if $q=7$, $9$ we have that $h\leq 2^{2g}$ and so, by theorem $\ref{madden}$, $h\leq 4$. $\blacksquare$

The above theorem tells us that for the function fields whose group of ideal classes has exponent two, it is enough to analyze only the cases when $h=2$, $4$, $8$, $16$, and $32$.

\section{Class number for $1\leq g \leq 10$}

In order to compute the class number $h$ of $K$, we describe the \textit{Artin-Weil zeta function}, but before that, we introduce some notation that will be used throughout this work. For a divisor $\mathcal{D}$, we put 
\[N(\mathcal{D})=q^{\mathrm{deg}\,\mathcal{D}}.\]

For a function field $K$ over the finite field of constants $\mathbb{F}_{q}$, $q=p^{n}$, the \textit{Artin-Weil zeta function} of $K/\mathbb{F}_{q}$ is defined, for $s\in \mathbb{C}$ with $\Re {s}>1$, by 
\[\zeta _{K}(s):=\sum_{\mathcal{D}}N(\mathcal{D})^{-s}=\prod_{\mathcal{P}}(1-N(\mathcal{P})^{-s})^{-1},\]
where the sum is taken over the \textit{positive} divisors $\mathcal{D}$, and the product is taken over all the places $\mathcal{P}$ of $K/\mathbb{F}_{q}$. The main results about this zeta function are contained in the following theorem.

\begin{thm}
The zeta function of $K/\mathbb{F}_{q}$ can be represented in the form 
\[\zeta _{K}(t)=\frac{L(t)}{(1-t)(1-qt)},\]
where $t=q^{-s}$, and $L(t)$ is a polynomial in $Z[t]$ of degree $2g$, $g$ being the genus of $K$. Furthermore $L(0)=1$ and $L(1)=h$ the class number of $K/\mathbb{F}_{q}$.
\end{thm}

The Riemann Hypothesis is the statement that all the zeros of $\zeta_K(s)$ are located in the line $\Re{s} = 1/2$. This is equivalent to the following theorem proved by Weil.

\begin{thm}[Riemann Hypothesis]
The reciprocal of the roots of $L(t)$ satisfy 
\[ |\alpha_{i}|=q^{1/2} \textrm{ for } i=1,\ldots ,2g \]
\end{thm}

Let 
\[L(t)=1+a_{1}t+a_{2}t^{2}+\cdots +a_{2g}t^{2g}=\prod_{i=1}^{2g}(1-\alpha_{i}t),\]
be the numerator of the zeta function. Then,
\begin{equation}
t^{-g}L(t)=t^{-g}+a_{1}t^{-(g-1)}+\cdots +a_{g}+\cdots +a_{2g}t^{g},\label{eq:1}
\end{equation}
is invariant if we replace $t$ by $q^{-1}t^{-1}$, i.e., $t^{-g}L(t)=q^{g}t^{g}L(\frac{1}{qt})$. Using this functional equation, we obtain 
\[a_{2g}=q^{g},a_{2g-1}=q^{g-1}a_{1},\ldots ,a_{g+1}=qa_{g-1},\]
and, hence it follows from $(\ref{eq:1})$ that
\begin{eqnarray}
L(t)=1+a_{1}t+a_{2}t^{2}+\cdots +a_{g}t^{g}+qa_{g-1}t^{g+1} \nonumber \\
  +\cdots +q^{g-1}a_{1}t^{2g-1}+q^{g}t^{2g}.\label{eq:2}
\end{eqnarray}
The expression for the class number is 
\begin{equation}
h=L(1)=(q^{g}+1)+a_{1}(q^{g-1}+1)+\cdots +a_{g-1}(q+1)+a_{g}.  \label{eq:cn}
\end{equation}
We shall calculate $a_{1}$, $a_{2}$, $a_{3}$, $a_{4}$, $a_{5}$, $a_{6}$, $a_{7}$, $a_{8}$, $a_{9}$ and $a_{10}$ which will suffice for the discussion of the cases $g=1$, $2$, $3$, $4$, $5$, $6$, $7$, $8$, $9$, $10$. Let $S_{\nu }=\sum_{i=1}^{2g}\alpha _{i}^{\nu }$. Then, 
\begin{equation}
-S_{\nu }=\sum_{d|\nu }d(N_{d}-n_{d}),  \label{eq:3}
\end{equation}
where $N_{d}$, $n_{d}$ denote respectively, the number of places of degree $d $ of $K$ and $k(x)$. The $a_{i}^{\prime }s$ are elementary symmetric functions of the $\alpha _{i}^{\prime }s$. In general, if $R$ is any commutative ring with identity, and $z$ is an indeterminate over the polynomial ring $R[x_{1}$, $\ldots $, $x_{n}]$, if $f(z)=(z-x_{1})(z-x_{2})\cdots (z-x_{n})$ then, 
\[ f(z)=z^{n}-\sigma _{1}z^{n-1}+\sigma _{2}z^{n-2}+\cdots +(-1)^{n}\sigma_{n},\]
with 
\[\sigma _{k}=\sigma _{k}(x_{1},\ldots ,x_{n})=\sum_{1\leq i_{1}<\cdots <i_{k}\leq n}x_{i_{1}}\cdots x_{i_{k}}\;\;\;(k=1,2,\ldots ,n).\]
The polynomial $\sigma _{k}$ is called the $k$\emph{th elementary symmetric polynomial} in the indeterminates $x_{1}$, $\ldots $, $x_{n}$ over $R$.

\begin{thm}[Newton's Recursion Formula]
Let $\sigma _{1}$, $\ldots $, $\sigma _{n}$ be the elementary symmetric polynomials in $x_{1}$, $\ldots $, $x_{n}$ over $R$, and let $S_{0}=n\in \mathbb{Z}$ and 
\[S_{k}=S_{k}(x_{1},\ldots ,x_{n})=x_{1}^{k}+\cdots +x_{n}^{k}\in R[x_{1},\ldots ,x_{n}],\]
for $k\geq 1$. Then the formula 
\[S_{k}-S_{k-1}\sigma _{1}+S_{k-2}\sigma _{2}+\cdots +(-1)^{m}\frac{m}{n}S_{k-m}\sigma _{m}=0\]
holds for $k\geq 1$, where $m=\mathit{min}(k$, $n)$.
\end{thm}

Using the Newton formula for $S_{\nu }$ in terms of $S_{1}$, $\ldots $, $S_{\nu -1}$ and the elementary symmetric functions, we obtain from $(\ref{eq:1})$:
\begin{eqnarray}
a_{1} &=&-S_{1},  \nonumber \\
a_{2} &=&\frac{-S_{2}+S_{1}^{2}}{2},  \nonumber \\
a_{3} &=&\frac{-2S_{3}+3S_{2}S_{1}-S_{1}^{3}}{6},  \nonumber \\
a_{4} &=&\frac{-6S_{4}+8S_{3}S_{1}+3S_{2}^{2}-6S_{2}S_{1}^{2}+S_{1}^{4}}{24},  \nonumber \\
a_{5} &=&\frac{-24S_{5}+30S_{4}S_{1}+20S_{3}S_{2}-20S_{3}S_{1}^{2}-15S_{2}^{2}S_{1}+10S_{2}S_{1}^{3}-S_{1}^{5}
}{120}.  \nonumber \\
\nonumber
\end{eqnarray}
In a similar way, we can find formulas for $a_{6},a_{7},a_{8},a_{9}$ and $a_{10}$.

Also, we have Dedekind's formulae, 
\begin{equation}
n_{d}=\left\{ 
\begin{array}{ll}
q+1 & \mbox{if $d=1$} \\[5mm] 
\frac{1}{d}\sum_{f|d}q^{f}\mu (\frac{d}{f}) & \mbox{if $d>1$}\label{df},
\end{array}
\right.
\end{equation}
where $\mu (m)$ denotes the M\"{o}bius function. It follows from $(\ref{eq:3})$ and $(\ref{df})$ that 
\begin{eqnarray*}
S_{1} &=&(q+1)-N_{1}, \\
S_{2} &=&(q^{2}+1)-N_{1}-2N_{2}, \\
S_{3} &=&(q^{3}+1)-N_{1}-3N_{3}, \\
S_{4} &=&(q^{4}+1)-N_{1}-2N_{2}-4N_{4}, \\
S_{5} &=&(q^{5}+1)-N_{1}-5N_{5}. \\ \nonumber
\end{eqnarray*}
In a similar way, we can find formulas for $S_{6},S_{7},S_{8},S_{9}$ and $S_{10}$.

The substitution of these values in the $a_{i}^{\prime }s$ gives, after simplification: 
\begin{eqnarray*}
a_{1} &=&N_{1}-q-1, \\
a_{2} &=&\frac{2N_{2}-N_{1}+N_{1}^{2}-2N_{1}q+2q}{2}, \\
a_{3} &=&\frac{%
-6N_{2}-3N_{1}^{2}q-6N_{2}q+N_{1}^{3}+6N_{2}N_{1}-N_{1}+6N_{3}+3N_{1}q}{6}%
, \\
a_{4} &=&\frac{%
12N_{2}-24N_{2}N_{1}q-N_{1}^{2}+24N_{2}q+N_{1}^{4}+12N_{2}^{2}+12N_{2}N_{1}^{2}-4N_{1}^{3}q%
}{24} \\
&&\frac{%
+24N_{3}N_{1}-24N_{3}q+24N_{4}+2N_{1}^{3}-12N_{2}N_{1}-2N_{1}-24N_{3}+4N_{1}q%
}{24}, \\
a_{5} &=&\frac{%
5N_{1}^{3}+5N_{1}^{2}q+40N_{2}N_{1}+60N_{2}N_{1}q-120N_{4}q+120N_{3}N_{2}+60N_{3}N_{1}^{2}%
}{120} \\
&&\frac{%
-60N_{2}^{2}q+20N_{2}N_{1}^{3}+60N_{2}^{2}N_{1}-5N_{1}^{4}q+120N_{4}N_{1}+N_{1}^{5}+5N_{1}^{4}%
}{120} \\
&&\frac{%
-60N_{2}^{2}-10N_{1}^{3}q-60N_{3}N_{1}+120N_{3}q-60N_{2}q-6N_{1}-60N_{2}%
}{120} \\
&&\frac{-120N_{4}+120N_{5}-120N_{3}N_{1}q-60N_{2}N_{1}^{2}q+10N_{1}q-5N_{1}^{2}}{120},
\\
\nonumber
\end{eqnarray*}
In a similar way, we can find formulas for $a_{6},a_{7},a_{8},a_{9}$ and $a_{10}$.

Substitution in $(\ref{eq:2})$, gives the numerator of the zeta function for 
\[ g=1,2,3,4,5,6,7,8,9,10. \]
Specifically, we obtain the following expressions for the class number.

$(g=1)$
\[h=N_{1}.\]

$(g=2)$
\[h=\frac{N_{1}-2q+2N_{2}+N_{1}^{2}}{2}.\]

$(g=3)$
\[h=\frac{3N_{1}^{2}-6N_{1}q+2N_{1}+6N_{3}+N_{1}^{3}+6N_{2}N_{1}}{6}.\]

$(g=4)$
\begin{eqnarray*}
h &=& \frac{ 6N_{1}+11N_{1}^{2}-12N_{1}q+12N_{2}+12N_{2}N_{1}-24N_{2}q-12N_{1}^{2}q}{24}  \nonumber \\
  && \frac{+24N_{4}+6N_{1}^{3}+24N_{3}N_{1}+12N_{2}^{2}+12N_{2}N_{1}^{2}+N_{1}^{4}}{24}.   \label{eqcn}
\end{eqnarray*}

$(g=5)$
\begin{eqnarray*}
h &=&\frac{-120N_{2}N_{1}q-40N_{1}q+24N_{1}+50N_{1}^{2}+100N_{2}N_{1}-60N_{1}^{2}q+35N_{1}^{3}}{120} \\
&&\frac{+60N_{1}N_{3}-120N_{3}q+60N_{2}N_{1}^{2}-20N_{1}^{3}q+10N_{1}^{4}+120N_{4}N_{1}}{120} \\
&&\frac{+120N_{3}N_{2}+60N_{3}N_{1}^{2}+60N_{2}^{2}N_{1}+20N_{2}N_{1}^{3}+120N_{5}+N_{1}^{5}}{120}.
\end{eqnarray*}

In a similar way, we can find formulas for $h$ when the genus is $6$, $7$, $8$, $9$, and $10$.

\section{Class number 2, 4, 8, 16, 32}

In this section we study function fields for which $h=2, 4, 8, 16, 32$. Using the Riemann Hypothesis, we have from equation $(\ref{eq:2})$ 
\begin{eqnarray}
L(t)&=&\prod_{j=1}^g(1-\sqrt{q}e^{i\theta_j}t)(1-\sqrt{q}e^{-i\theta_j}t) \nonumber \\
&=&\prod_{j=1}^g(1-2\sqrt{q}t\cos{\theta_j}+qt^2).  \label{lhr} 
\end{eqnarray}
Since $L(1)=h$, we obtain from $(\ref{lhr})$ 
\begin{eqnarray}
h \geq (\sqrt{q}-1)^{2g}.  \label{dh}
\end{eqnarray}

\begin{rem}
\label{ccn1} $K$ can not have a number of places of degree one bigger than $h$, unless its genus is zero. For the existence of more than $h$ places of degree one, would imply that two of them, say $\mathcal{P}_{1}$ and $\mathcal{P}_{2}$, belong to the same class, i.e., $\mathcal{P}_{1}-\mathcal{P}_{2}=(x)$ and $[K:k(x)]=\mathrm{deg}\,(x)_{0}=1$. Also, fields of genus zero have class number one. Therefore, $N_{1}\leq h$.
\end{rem}

\subsection{Class number 2}

It follows from $(\ref{dh})$ that $h > 2$ if $q \geq 6$. Hence, $q=2, 3, 4, 5$. Then, it is proved in~\cite{cn2}

\begin{thm}
\label{h2}Let $K/k$ be a function field of genus $g$, with $k=\mathbb{F}_{q}$ the field of constants. Then, its class number is larger than two if any of the following conditions are satisfied.\newline
$\left( i\right) $ $q=4$ or $5$ and $g\geq 2$.\newline
$\left( ii\right) $ $q=3$, $g\geq 3$.\newline
$\left( iii\right) $ $q=2$, $g\geq 6$.
\end{thm}

\subsection{Class number 4}

We proceed in the same form that for the case $h=2$, in the way that did it~\cite{cn2}. It follows from $(\ref{dh})$ that $h>4$ if $q\geq 10$. Hence, $q=2,3,4,5,7,8,9$. Then, we have

\begin{thm}
\label{h4}Let $K/k$ be a function field of genus $g$, with $k=\mathbb{F}_{q}$ the field of constants. Then, its class number is larger than four if any of the following conditions are satisfied.\newline
$\left( i\right) $ $q=5$ or $7,8,9$ and $g\geq 2$.\newline
$\left( ii\right) $ $q=4$, $g\geq 3$.\newline
$\left( iii\right) $ $q=3$, $g\geq 4$.\newline
$\left( iv\right) $ $q=2$, $g\geq 7$.
\end{thm}

{\bf Proof. }
We follow the proof in~\cite{cn2} step by step. The Riemann Hypothesis is equivalent to the inequality 
\begin{equation}
|N_{1}-(q+1)|\leq 2g\sqrt{q}.  \label{erh}
\end{equation}
Let $\overline{K}/\overline{k}$ be a constant field extension of $K/k$ of degree $2g-1$. Then, $\overline{K}/\overline{k}$ has the same genus as $K/k$ . We apply $(\ref{erh})$ to $\overline{K}/\overline{k}$ and obtain 
\begin{equation}
\overline{N}_{1}\geq q^{2g-1}+1-2gq^{(2g-1)/2}.  \label{erho}
\end{equation}
A place of degree $d$ of $K$ decomposes (\cite{aff}, p. 163) in $\overline{K}$ as $(d,2g-1)$ places of degree $\frac{d}{(d,2g-1)}$. Places of degree one of $\overline{K}$, therefore, lie over such places of $K$ of which degree divides $2g-1$. Using $(\ref{erho})$, it follows that $K$ has, at least 
\[ \frac{q^{2g-1}+1-2gq^{(2g-1)/2}}{2g-1}\]
positive divisors of degree $2g-1$: If $N=$ positive divisors of degree $2g-1 $, then 
\[N\geq \sum_{d\mid 2g-1}N_{d},\]
and therefore 
\begin{eqnarray*}
\left( 2g-1\right) N &\geq &\sum_{d\mid 2g-1}\left( 2g-1\right) N_{d} \\
&\geq &\sum_{d\mid 2g-1}dN_{d}=\overline{N}_{1} \\
&\geq &q^{2g-1}+1-2gq^{(2g-1)/2}.
\end{eqnarray*}
On the other hand, the total number of positive divisors in $K$ of degree $2g-1$ is $\frac{h}{q-1}(q^{g}-1)$ (\cite{aff}, p. 160). Therefore, $h$ is larger than four if 
\[(q-1)(q^{2g-1}+1-2gq^{(2g-1)/2})>4(q^{g}-1)(2g-1).\]
We see that $h>4$ if 
\[S(q,g)=(q-1)(q^{2g-1}+1-2gq^{(2g-1)/2})-4(q^{g}-1)(2g-1)>0.\]
This occurs in the cases listed in the theorem. $\blacksquare$

\subsection{Class number 8}

It follows from $(\ref{dh})$ that $h > 8$ if $q \geq 15$. Hence, $q=2, 3, 4, 5,7, 8, 9, 11, 13$. First, we have

\begin{thm}
\label{h8}Let $K/k$ be a function field of genus $g$, with $k=\mathbb{F}_{q}$ the field of constants. Then, its class number is larger than eight if any of the following conditions are satisfied.\newline
$\left( i\right) $ $q=7$ or $8,9,11,13$ and $g\geq 2$.\newline
$\left( ii\right) $ $q=5$, $g\geq 3$.\newline
$\left( iii\right) $ $q=4$, $g\geq 4$.\newline
$\left( iv\right) $ $q=3$, $g\geq 5$.\newline
$\left( v\right) q=2$, $g\geq 9$.
\end{thm}

{\bf Proof. }
Same proof as in Theorem $\ref{h4}$ with 
\[S(q,g)=(q-1)(q^{2g-1}+1-2gq^{(2g-1)/2})-8(q^{g}-1)(2g-1). \blacksquare \]

\subsection{Class number 16}

It follows from $(\ref{dh})$ that $h>16$ if $q\geq 26$. Hence,
\[ q=2,3,4,5,7,8,9,11,13,16,17,19,23,25.\]
 First, we have

\begin{thm}
\label{h16}Let $K/k$ a function field of genus $g$, with $k=\mathbb{F}_{q}$ the field of constants. Then, its class number is larger than sixteen if any of the following conditions are satisfied.\newline
$\left( i\right) $ $q=9$ or $11,13,16,17,19,23,25$ and $g\geq 2$.\newline
$\left( ii\right) $ $q=5$ or $7,8$ and $g\geq 3$.\newline
$\left( iii\right) $ $q=4$, $g\geq 4$.\newline
$\left( iv\right) $ $q=3$, $g\geq 5$.\newline
$\left( v\right) $ $q=2$, $g\geq 10$.
\end{thm}

{\bf Proof. }
Same proof as in Theorem $\ref{h4}$ with 
\[ S(q,g)=(q-1)(q^{2g-1}+1-2gq^{(2g-1)/2})-16(q^{g}-1)(2g-1). \blacksquare  \]

\subsection{Class number 32}

It follows from $(\ref{dh})$ that $h>32$ if $q\geq 45$. Hence, 
\[q=2,3,4,5,7,8,9,11,13,16,17,19,23,25,27,29,31,32,37,41,43.\]
First, we have

\begin{thm}
\label{h32}Let $K/k$ a function field of genus $g$, with $k=\mathbb{F}_{q}$ the field of constants. Then, its class number is larger than thirty two if any of the following conditions are satisfied.\newline
$\left( i\right) $ $q=11$ or $13,16,17,19,23,25,27,29,31,32,37,41,43$ and $g\geq 2$.\newline
$\left( ii\right) $ $q=7$ or $8,9$ and $g\geq 3$.\newline
$\left( iii\right) $ $q=5$, $g\geq 4$.\newline
$\left( iv\right) $ $q=4$, $g\geq 5$.\newline
$\left( v\right) $ $q=3$, $g\geq 6$.\newline
$\left( vi\right) $ $q=2$, $g\geq 11$.
\end{thm}

{\bf Proof. }
Same proof as in Theorem \ref{h4} with 
\[S(q,g)=(q-1)(q^{2g-1}+1-2gq^{(2g-1)/2})-32(q^{g}-1)(2g-1). \blacksquare  \]

\subsection{Initial modification to Madden's theorem}

From the proof of the theorem $\ref{menorque}$ and theorems $\ref{h2}$, $\ref {h4}$, $\ref{h8}$, $\ref{h16}$ and $\ref{h32}$, we can give an initial modification to the theorem $\ref{madden}$.

\begin{thm}[Madden, initial modification]
If $K$ is a qua- \newline dratic extension of $k(x)$ of genus $g$, where $k$ is a finite field of order $q$, in which the infinite place ramifies, and if the ideal class group of $K$ has exponent two, then:\newline
(a) h=2.\newline
\hspace*{.3in}(i) $q=4$, $5$, g=1. \newline
\hspace*{.3in}(ii) q=3, $1\leq g \leq 2.$ \newline
\hspace*{.3in}(iii) q=2, $1\leq g \leq 5.$ \newline
(b) h=4.\newline
\hspace*{.3in}(i) $q=5$, $7$, $9$, g=1. \newline
\hspace*{.3in}(ii) q=4, $1\leq g \leq 2.$ \newline
\hspace*{.3in}(iii) q=3, $1\leq g \leq 3.$ \newline
\hspace*{.3in}(iv) q=2, $2\leq g \leq 6.$ \newline
(c) h=8.\newline
\hspace*{.3in}(i) q=5, $1\leq g \leq 2.$ \newline
\hspace*{.3in}(ii) q=3, $1\leq g \leq 4.$ \newline
\hspace*{.3in}(iii) q=2, $4\leq g \leq 8.$ \newline
(d)  h=16.\newline
\hspace*{.3in}(i) q=5, $1\leq g \leq 2.$ \newline
\hspace*{.3in}(ii) q=3, $1\leq g \leq 4.$ \newline
\hspace*{.3in}(iii) q=2, $7\leq g\leq 8.$ \newline
(e) h=32.\newline
\hspace*{.3in}(i) q=3, $1\leq g \leq 4.$ \newline
\end{thm}

{\bf Proof. }
Let $q=2$, then $h\leq 16$. If $h=2$, we have that some place different from the infinite place must ramify; from here, by the equation $(\ref{mayor 0})$ we have that $2g+1>2\left( 1\right) $, and so, $g\geq 1$. If $h=4$, then must ramify two places differents from the infinite place; hence, by equation $(\ref{mayor 0})$ we have that $2g+1>2\left( 1+1\right) $, therefore, $g\geq 2$. If $h=8$, then must ramify three places; again by equation $(\ref {mayor 0})$ we have that $2g+1>2\left( 1+1+2\right) $, so, $g\geq 4$. If $h=16 $, four places must ramify, so $2g+1>2\left( 1+1+2+3\right) $, and therefore, $g\geq 7$.
The other results are similarly obtained. $\blacksquare$

\section{Quadratic function fields with exponent two ideal class group}

In this part we give the full list of examples of quadratic function fields with exponent two ideal class group.

\subsection{Relation between places of degree r and places of degree one}

Using the Hasse-Weil bound $(\ref{erh})$ it is possible to give an estimation for the number of places of fix degree $r$.

Given a function field $K/\mathbb{F}_{q}$ of genus $g$, denote 
\[N_{r}=\left| \left\{ P\in \mathbb{P}_{F}:\deg P=r\right\} \right| .\]
There is a close relation between the numbers \thinspace $N_{r}$ and $nps$ (the number of places of degree one in the extension $K_{s}=K\mathbb{F}_{q^{s}}$) given by: 
\[ npr=\sum_{d\mid r}d\cdot N_{d} \]
(the sum is over all the divisors $d\geq 1$ that divide $r$).

The M\"{o}bius inversion formula transforms the last equation in: 
\begin{equation}
r\cdot N_{r}=\sum_{d\mid r}\mu \left( \frac{r}{d}\right) \cdot npd.  \label{Nr}
\end{equation}

\subsection{Discussion of an example when $q=2$ and $g=2$}
We know that $h\leq 2^{g}$; in this case $h\leq 2^{2}=4$. So, $h=2$ or $4$. If $h=2$, by theorem $\ref{caracterizacion}$, $t-1=1$ and hence, $t=2$. If $h=4$, by theorem $\ref{caracterizacion}$, $t-1=2$ and $t=3.$ Therefore, using the Hasse Normal Form,
\[u=\frac{q\left( x\right) }{x^{\gamma _{1}}\left( x+1\right) ^{\gamma_{2}}\left( x^{2}+x+1\right) ^{\gamma _{3}}}\]
and 
\[\deg q\left( x\right) =5-\sum \deg p_{i}\left( x\right).\]
First suppose that $h=2$. There are several possibilities for $u$.
We discuss only the case when the Hasse Normal Form for $u$ is: 
\[u=\frac{q\left( x\right) }{x^{\gamma _{1}}} \]
and 
\[\deg q\left( x\right) =4\textmd { with }q\left( 0\right) \neq 0.\]
Since $y^{2}+y=u$ we have 
\begin{eqnarray*}
y^{2}+xy &=&x\cdot q\left( x\right) , \\
y^{2}+x^{2}y &=&x\cdot q\left( x\right) .
\end{eqnarray*}
Next we show how to find all the function fields of the form $y^{2}+xy=x\cdot q\left( x\right) :$ Note that $q\left( x\right) $ is of the form 
\[q\left( x\right) =x^{4}+a_{3}x^{3}+a_{2}x^{2}+a_{1}x+1.\]
Now we give values to $\left( a_{3,}a_{2},a_{1}\right) $ over $\mathbb{F}_{2}^{3}$, and count the places of degree one and two using $(\ref{Nr})$. If $h=2$, where $h$ is given by $(\ref{eqcn})$, then we have a legitimate example of a function field with exponent two ideal class group.

\subsection{Discussion of an example when $q=3$ and $g=1$}

This is an example of an elliptic function field. By theorem $\ref{nec eff}$, $y^{2}=f\left( x\right) $ where $f\left( x\right) $ is a square free polynomial of degree three. There are two possibilities: $f\left( x\right)$ factorizes in three polynomials of degree one or one polynomial of degree one and the other of degree two. In the first case, there are a four places that ramify and so $t=4$. By theorem $\ref{caracterizacion}$, we are looking for examples where $h=4$. In this case, our general curve is $y^{2}=a\left(x\right) \left( x+1\right) \left( x+2\right) $ with $a\in \mathbb{F}_{3}^{\times }$. Varying $a$ we find the curves that have $h=4$. In the second case, there are three places that ramify and so $t=3$. By theorem $\ref{caracterizacion}$, we are looking for examples where $h=2$. In this case, our general curve is $y^{2}=a\left( x+b\right) \left( x+cx+d\right)$ with $a$, $b$, $c$, $d\in \mathbb{F}_{3}^{\times }$.

\subsection{Discussion of an example when $q=3$ and $g=2$}

This is an example of an hyperelliptic function field. By theorem $\ref {hyper eff}$, $y^{2}=f\left( x\right) $ where $f\left( x\right) $ is a square free polynomial of degree five or six. The case in which the degree of $f\left( x\right) $ is six, does not occur because in this situation $\infty $ does not ramify. So the degree of $f\left( x\right) $ is five. For $h=2$, according to \cite{cn2}, there are no examples. If we want $h=4=2^{2}$, by theorem $\ref{caracterizacion}$, $h=2^{t-2}$ and from here, $t=4$. Since $\infty $ must ramify, $t=3$. So $f\left( x\right) $ must decompose as a product of three polynomials. Using degrees $1+1+3$ or $1+2+2$ polynomials, the only possibilities for $f\left( x\right) $ are:  
\begin{eqnarray*}
y^{2} &=&a\left( x\right) \left( x+b\right) \left( x^{3}+cx^{2}+dx+e\right) ,
\\
y^{2} &=&a\left( x+1\right) \left( x+2\right) \left(
x^{3}+bx^{2}+cx+d\right) , \\
y^{2} &=&a\left( x+b\right) \left( x^{2}+bx+c\right) \left(
dx^{2}+ex+f\right).
\end{eqnarray*}
For $h=8$, the only possibility for $f\left( x\right) $ is $a\left( x\right) \left( x+1\right) \left( x+2\right) \left( x^{2}+bx+c\right) $. Note that the cases $h=16$ and $h=32$ can not happen.

\subsection{Full list of examples up to $\mathbb{F}_{q}$-isomorphism}

\subsubsection{Class number 2}

\begin{exmp}$q=2$, $g=1$, $N_{1}=2$ \\ 
\begin{itemize}
\item  $y^{2}+xy+x\left( x^{2}+x+1\right) =0.$
\end{itemize}
\end{exmp}

\begin{exmp}$q=2$, $g=2$, $N_{1}=2$, $N_{2}=1$ \\
\begin{itemize}
\item  $y^{2}+xy+x\left( x^{4}+x+1\right) =0.$
\end{itemize}
\end{exmp}

\begin{exmp}$q=2$, $g=2$, $N_{1}=1$, $N_{2}=3$ \\
\begin{itemize}
\item  $y^{2}+(x^{2}+x+1)y+(x^{2}+x+1)(x^{3}+x+1)=0.$
\end{itemize}
\end{exmp}

\begin{exmp}$q=2$, $g=3$, $N_{1}=1$, $N_{2}=3$, $N_{3}=0$ \\
\begin{itemize}
\item  $y^{2}+(x^{2}+x+1)y+(x^{2}+x+1)(x^{5}+x^{2}+1)=0.$
\end{itemize}
\end{exmp}

\begin{exmp}$q=2$, $g=3$, $N_{1}=1$, $N_{2}=2$, $N_{3}=1$ \\ 
\begin{itemize}
\item  $y^{2}+(x^{3}+x^{2}+1)y+(x^{3}+x^{2}+1)(x^{4}+x^{3}+1)=0.$
\end{itemize}
\end{exmp}

\begin{rem}
The Example $A$ of $\left[ \textmd {LMQ75, page 25}\right] $, 
\[
y^{2}+\left( x^{3}+x+1\right) y=\left( x^{3}+x+1\right) \left(
x^{4}+x+1\right) 
\]
has $g=3$ and they claim that $h=2$ (with $N_{1}=1$, $N_{2}=2$, $N_{3}=1$). But our computations show that this example has $h=4$ (with $N_{1}=1$, $N_{2}=2$, $N_{3}=3$). The correct example is 
\[y^{2}+\left( x^{3}+x+1\right) y=\left( x^{3}+x+1\right) \left(x^{4}+x^{3}+x^{2}+x+1\right)\]
and here, $h=2$ (with $N_{1}=1$, $N_{2}=2$, $N_{3}=1$).
\end{rem}

\begin{exmp} $q=3$, $g=1$, $N_{1}=2$ 
\begin{itemize}
\item  $y^{2}-(x+2)(x^{2}+1)=0.$
\end{itemize}
\end{exmp}

\begin{exmp}$q=4$, $g=1$, $N_{1}=2$\\
Let $\alpha $ be a root of the polynomial $x^{2}+x+1$.
\begin{itemize}
\item  $y^{2}+yx+\alpha x(x^{2}+x\alpha +\alpha )=0.$
\end{itemize}
\end{exmp}

\begin{exmp}$q=5,g=1,N_{1}=2$ \\
\begin{itemize}
\item  $y^{2}-x(x^{2}+2)=0.$
\end{itemize}
\end{exmp}

\subsubsection{Class number 4}

\begin{exmp}$q=2,g=2,N_{1}=3,N_{2}=0$ \\ 
\begin{itemize}
\item  $y^{2}+x(x+1)y+x(x+1)(x^{3}+x+1)=0.$
\end{itemize}
\end{exmp}

\begin{exmp}$q=2,g=3,N_{1}=3,N_{2}=0,N_{3}=0$ \\ 
\begin{itemize}
\item  $y^{2}+x(x+1)y+x(x+1)(x^{5}+x^{3}+x^{2}+x+1)=0.$
\end{itemize}
\end{exmp}

\begin{exmp}$q=2,g=3,N_{1}=2,N_{2}=2,N_{3}=0$ \\
\begin{itemize}
\item  $y^{2}+x(x^{2}+x+1)y+x(x^{2}+x+1)(x^{4}+x+1)=0.$
\end{itemize}
\end{exmp}

\begin{exmp}$q=3,g=2,N_{1}=3,N_{2}=1$ \\
\begin{itemize}
\item  $y^{2}-2x(x+1)(x^{3}+x^{2}+x+2)=0.$
\end{itemize}
\end{exmp}

\begin{exmp} $q=3,g=2,N_{1}=2,N_{2}=4$ \\
\begin{itemize}
\item  $y^{2}-2(x+2)(x^{2}+1)(x^{2}+x+2)=0.$
\end{itemize}
\end{exmp}

\subsubsection{Class number 8}

\begin{exmp}$q=3,g=2,N_{1}=4,N_{2}=1$ \\ 
\begin{itemize}
\item  $y^{2}-x(x+1)(x+2)(x^{2}+1)=0.$
\end{itemize}
\end{exmp}

\subsubsection{Class number 16}

Computer search shows that there is no such example.

\subsubsection{Class number 32}

Computer search shows that there is no such example.

\section{Final modification to Madden's theorem}

In the last section, we found the complete list of the fields whose ideal class group has exponent two. Our calculations suggest what could be the final form of Madden's theorem.

\begin{thm}[Madden, final modification]
If $K$ is a quadratic extension of $k(x)$ of genus $g$, where $k$ is a finite field of order $q$, in which the infinite place ramifies, and if the ideal class group of $K$ has exponent two, then:\newline
$\left( a\right) $ $h=2.\newline
\begin{array}{cc}
& \left( i\right) 
\end{array}
q=3,4,5$, $g=1.\newline
\begin{array}{cc}
& \left( ii\right) 
\end{array}
q=2$, $1\leq g \leq 3.\newline
\left( b\right) $ $h=4.\newline
\begin{array}{cc}
& \left( i\right) 
\end{array}
q=3$, $g=2.\newline
\begin{array}{cc}
& \left( ii\right) 
\end{array}
q=2, 2\leq g\leq 3.\newline
\left( c\right) $ $h=8.$\newline
$
\begin{array}{cc}
& \left( i\right) 
\end{array}
q=3$, $g=2.$\hfill $\square $
\end{thm}

\section*{Acknowledgements}
We are grateful to Dinesh Thakur for his comments and suggestions regarding the content and draft of this article.

\end{document}